\documentclass[a4paper]{article}

\usepackage{amsmath}
\usepackage{amsfonts}
\usepackage{amssymb}         
\usepackage{amsopn}
\usepackage{theorem}          
\usepackage{latexsym}
\usepackage{graphicx}        
\usepackage{mathrsfs}		
\usepackage{psfrag}
\usepackage{algorithmic}
\usepackage{algorithm}
\usepackage{color}
\usepackage{verbatim}
\usepackage{url}

\definecolor{offwhite}{gray}{0.92}
\definecolor{reasonablegreen}{rgb}{0,0.5,0}

\newtheorem{result}{\ }[section]
\theoremstyle{changebreak}                
\newtheorem{thm}[result]{Theorem}

\newtheorem{lem}[result]{Lemma}
\newtheorem{cor}[result]{Corollary}

 {{\nopagebreak\hspace*{\fill}$\Box$\par\vspace{12pt}}}
\newcommand{\transpose}[1]{{#1}^{\top}}

\newcommand{\supp}[1]{\mathsf{supp}({#1})}
\newcommand{\prob}[1]{\mathsf{Prob}\big({#1}\big)}
\newcommand{\Prob}[1]{\mathsf{Prob}\left[{#1}\right]}
\newcommand{\red}[1]{{\color{black}#1}}
\begin{document}

\thispagestyle{empty}
\begin{center} 

{\LARGE Decoding noisy messages: a method that just shouldn't work}
\par \bigskip
{\sc Leo Liberti${}^{1}$} 
\par \bigskip
\begin{minipage}{15cm}
\begin{flushleft}
{\small
\begin{itemize}
\item[${}^1$] {\it LIX CNRS, \'Ecole Polytechnique, Institut Polytechnique de Paris, F-91128 Palaiseau, France} \\ Email:\url{liberti@lix.polytechnique.fr}
\end{itemize}
}
\end{flushleft}
\end{minipage}
\par \medskip \today
\end{center}
\par \bigskip

\begin{abstract} 
  This paper is about receiving text messages through a noisy and costly line. Because the line is noisy we need redundancy, but because it is costly we can afford very little of it. I start by using well-known machinery for decoding noisy messages (compressed sensing), then I attempt to reduce the redundancy (using random projections), until I get to a point where I use more orthogonal vectors than the space dimension allows. Instead of grinding to a halt or spurting out noise, this method is still able to decode messages correctly or almost correctly. I have no idea why the method works: this is my first reason for writing this paper using a narrative instead of formal scientific style (the second one is that I am tired of writing semi-formal prose, and long for a change).
\end{abstract}


\section{The risks of office gossip}
\label{s:story}
At a time when I was mainly working on optimization and Euclidean distance geometry, I happened to walk past an open office door, where I heard someone mumble ``random projection''; my eye glanced upon a whiteboard bearing the formula
\begin{equation}
  \|Tx-Ty\|_2\approx\|x-y\|_2. \label{jll1}
\end{equation}
Intrigued by the Euclidean distances, I stopped, entered the office, and interrupted the ongoing discussion. ``What's that $T$?'', I asked. I was told that it was a $k\times n$ random matrix, and that the vectors $x,y$ were any pair of $n$-dimensional vectors out of a finite set. I immediately objected that no such $T$ could exist, unless the ``approximation'' relationship between the two distances was taken very liberally. ``Don't you know the Johnson-Lindenstrauss lemma?'', was the reply. I did not. After being given a summary explanation I left that office, my head in turmoil, and headed towards the library in search of more information.

Upon coming back from the library, I met a friendly colleague and asked him what he knew about random matrices. He immediately started talking about things unrelated to Eq.~\eqref{jll1}. He had a disposition to stray far from the subject of conversation, which was one of the reasons I liked him. This attitude resonated with my habit of making connections between different topics, the more disconnected the better. He was talking about tiny cameras with a ridiculously small number of pixels, but able to take fantastically detailed pictures nonetheless. He asked me ``do you know how many pixels I really need for those pictures? guess!''. I guessed in the thousands, then in the hundreds, then gave up. ``One!'', he said, smiling broadly. For the second time that day, I was nonplussed, and couldn't believe what I was hearing. A whole parallel world existed in my colleagues' minds, something that appeared to go against the expected tenets of normal mathematics, and I was cut off. This had to stop: I would measure the extent of their lies, and ascertain what truth might be hidden in their unlikely stories. 

Thus started a long period of learning, when I became aware of a large part of mathematics I had ignored so far. I initially focused on the office gossip about Euclidean distances: I read the lemma by Johnson and Lindenstrauss \cite{jllemma} (whose proof remained unintelligible to me for a long time), learned some unexpected properties of random matrices having components sampled from Gaussian distributions \cite{vempala06}, and became aware of their relation with data science \cite{dirksen}. I heard a reference to the one-pixel nonsense in a talk by Emmanuel Cand\`es at the International College of Mathematicians in 2014 \cite{candes}, where, again, random matrices were mentioned. A few web searches later, I had uncovered a name, that of Matou\v{s}ek, connected to these concepts.

Connections notwithstanding, I was by then fully convinced that my ignorance was abysmal, all-encompassing, and was procuring me a scientific damage beyond repair. I had been living a quiet mathematical life, sheltered from the stormy developments of the youngest greatgrandchildren of functional analysis, linear algebra, and geometry. Contemporary mathematics was passing me by. The Mathematical Programming (MP) niche that looked so rich and endless appeared like a prehistorical cave with a few ancient rat bones scattered about. 

Luckily, as I walked into one of my favorite maths bookshops (now unfortunately closed forever), I stumbled upon a copy of a book about Linear Programming (LP) written by Ji\v{r}\'{\i} Matou\v{s}ek and Bernd G\"artner \cite{matousek-lp}. It attracted my attention because it had an in-depth discussion of the Delsarte bound for codes \cite{delsarte}: my interest was in understanding some difficult papers about the Delsarte bound for \textit{spherical} codes \cite{delsarte2}, which arise in the study of the Kissing Number Problem \cite{knp_airo17}. Unfortunately, Matou\v{s}ek and G\"artner only discussed the application to discrete codes. The next section \cite[\S 8.5]{matousek-lp}, however, was titled ``Sparse solutions of linear systems'', which resonated with the title of Cand\`es' talk at ICM14, ``The mathematics of sparsity''. So, charmed by serendipity, I bought the book and kept on reading.

\section{Noisy communication channels}
\label{s:noisy}
At this point I will hand over the story to the authors of \cite{matousek-lp}.
\begin{quote}
  We begin by discussing error-correcting codes again, but this time we want to send a sequence $w\in\mathbb{R}^d$ of $d$ real numbers. Or rather not we, but a deep-space probe which needs to transmit its priceless measurements represented by $w$ back to Earth. We want to make sure that all components of $w$ can be recovered correctly even if some fraction, say 8\%, of the transmitted numbers are corrupted, due to random errors or even maliciously (imagine that the secret \textit{Brotherhood for Promoting the Only Truth} can somehow tamper with the signal slightly in order to document the presence of supernatural phenomena in outer space). We admit \textit{gross errors}; that is, if the number $3.1415$ is sent and it gets corrupted, it can be received as $2152.66$, or $3.1425$, or $-1011$, or any other real number.
  
  Here is a way of encoding $w$: We choose a suitable number $n > d$ and a suitable $n\times d$ encoding matrix $Q$ of rank $d$, and we send the vector $z = Qw \in\mathbb{R}^n$. Because of the errors, the received vector is not $z$ but $\bar{z} = z + \red{x}$, where \red{$x\in\mathbb{R}^n$} is a vector with at most $r = \lfloor 0.08n\rfloor$ nonzero components. We ask, under what conditions on $Q$ can $z$ be recovered from $\bar{z}$?
  
Somewhat counterintuitively, we will concentrate on the task of finding \red{the ``error vector'' \red{$x$}}. Indeed, once we know $x$, we can compute $w$ by solving the system of linear equations $Qw = z = \bar{z} - \red{x}$. The solution, if one exists, is unique, since we assume that $Q$ has rank $d$ and hence the mapping $w \mapsto Qw$ is injective.
\end{quote}

Let me summarize: we need to send a vector $w\in\mathbb{R}^d$ on a noisy channel. We would like to find an encoding $n\times d$ matrix $Q$, with $n>d$, and send $z=Qw\in\mathbb{R}^n$ to the recipient. We assume that both parties know $Q$. There is a low probability $\Delta$ of communication error (assumed to be 8\% in \cite{matousek-lp}), which implies that, on average, on sufficiently long vectors, $\Delta n$ components of $z$ will be wrong at reception. Let \red{$x\in\mathbb{R}^n$} be the error vector, so that the received message is $\bar{z}=z+\red{x}$. Evidently, $x$ will have on average only $\Delta n$ nonzero components. What conditions on $Q$ should we impose so that the recipient can retrieve $w$ from $\bar{z}$?

The book \cite{matousek-lp} goes on to explain that we (stepping in the recipient's shoes) should choose an $m\times n$ matrix $A$ with $m=n-d$ such that $AQ=0$, and let $b=A\bar{z}$. We should then find the sparsest solution $x'$ to the linear system $Ax=b$. Armed with $x'$, we find $z'=\bar{z}-x'$, and finally $w'=(\transpose{Q}Q)^{-1}\transpose{Q}z'$. In general, we should find that $w=w'$, or at least that $w,w'$ are not very different. \red{Note that the sparsity of a vector $x'\in\mathbb{R}^n$ is related to its support, i.e.~the set of indices of the nonzero components $\supp{x'}=\{j\le n\;|\;x'_j\not=0\}$: $x'$ is sparse when its support is small w.r.t.~$n$.}

There are three technical points worth making about the preceding discussion.
\begin{enumerate}
\item We compute $b$ as $A\bar{z}$ and then want to solve the system $Ax=b$. An encoded noisy message $\bar{z}$ and an error vector $x$ (two vectors with apparently very different meanings) are both being cast in the same role, in the sense that they are being multiplied by $A$ and give $b$ as a result: why? The reason is that we have
  \begin{equation}
    b = A\bar{z}=A(z+x)=A(Qw+x)=(AQ)w + Ax = Ax,\label{decode}
  \end{equation}
since $AQ=0$ by definition.
\item The system $Ax=b$ is underdetermined, since $A$ is an $(n-d)\times n$ matrix. It has a solution as long as the rank of $A$ is equal to the rank of the extended matrix $(A,b)$, a sufficient condition being that $A$ has full rank $n-d$.
\item We compute an approximation $w'$ of the original message as $(\transpose{Q}Q)^{-1}\transpose{Q}z'$: this is the well-known Moore-Penrose pseudoinverse, which is being used here to ``invert'' the (overdetermined, hence in general uninvertible) linear system $Qw=z$ used during the encoding.
\end{enumerate}
\red{Note that, in the aforementioned discussion, $x$ is being used both as a numerical error vector in $\mathbb{R}^n$ and as symbol for a vector of variables in the system $Ax=b$ to be solved. To avoid ambiguities, I shall henceforth reserve $x$ for the variable vector (as well as to make occasional appearances as a vector in theorem statements), and use $\bar{x}$ for the error vector. As long as I am discussing notation, let me clarify that $|\cdot|$ denotes both absolute value and set cardinality: I accept this ambiguity because both notations are standard and established, and because I made it clear what sign means what in the only possible case of ambiguity, Eq.~\eqref{eq:amb}.}

And then, of course, there is the metaphorical elephant in the proverbial room: why are we finding the sparsest solution $x'$ to $Ax=b$, and why, having found it, would it be equal to the error vector $\bar{x}$, which does not seem likely that we --- the receiving party, now --- could ever know, since it depends on the whims of the noisy line?

\section{The one-pixel camera}
\label{s:1norm}
Eventually, I found some answers to the last questions of the foregoing section. But they are not the type of answers I had come to expect from mathematics. Rather than the short, clear-cut first-order logical propositions with the usual, well-known, and trusted sentences of the \red{Zermelo-Fraenkel-Choice (ZFC)} axiom system \cite{kunen}, which are either provable or unprovable (and those that are provable end up being true in all sensible models), I found myself entangled in vague probabilistic statements. Their associated probabilities exceed certain threshold values that can be made arbitrarily close to one in function of a monotonically varying parameter, the actual value of which is at best hard to compute, and at worst unknowable, in all practical settings. I found something similar in the study of random projections I had pursued a couple of years earlier --- something else that Johnson, Lindenstrauss and the one-pixel camera had in common.

Let us start from one of the basic questions: how does one find the sparsest solution to an underdetermined linear system $Ax=b$? A MP formulation for this problem is $\min\{\|x\|_0\;|\;Ax=b\}$, where $\|\cdot\|_0$ is the ``zero-norm'', which is not really a norm, but simply counts the number of nonzero components of $x$. For later reference, I shall call this formulation $P^0(A,b)$. It turns out that $P^0(A,b)$ is \textbf{NP}-hard by reduction from \textsc{Exact Cover by 3-Sets} \cite[A6{[MP5]}]{gareyjohnson}. 

The signal processing community had been trying to solve these noisy decoding problems for decades: for example in order to build telephones, radios, TVs. In attempting to solve $P^0(A,b)$ efficiently they tried everything, including the replacement of the zero-norm with the much friendlier $\ell_1$ norm. The corresponding problem $\min\{\|x\|_1\;|\;Ax=b\}$, which I will refer to as $P^1(A,b)$, can be reformulated exactly to the following LP:
\begin{equation}
  \left.\begin{array}{rrcl}
    \min & \sum\limits_{j\le n} & s_j & \\
    \forall j\le n & -s_j \le & x_j &\le s_j \\
    & Ax &=& b.
  \end{array}\right\} \label{basispursuit}
\end{equation}
With much amazement, signal processing researchers noticed that, if they chose the dimensions $d$ and $n$ carefully, the solution \red{$x'$ of $P^1(A,b)$} was the same as the solution of $P^0(A,b)$. This happened too often for them to accept it as a ``coincidence''. The theoretical justification came from the work of many people, but principally David Donoho \cite{donoho}, Terence Tao, and Emmanuel Cand\`es \cite{candestao2005,candes_tao}. I do not really want to provide a full proof here: I will just sketch a proof strategy which I borrowed from many sources, but mainly from \cite{damelin}.

Part of the vagueness I referred to above comes from a practical need: that of computing with floating point numbers. If you find that a component of \red{$x'$ is $O(10^{-9})$}, do you say that it is zero or nonzero? Let us provide ourselves with a small $\epsilon>0$. We say that a scalar $\sigma$ is ``almost zero'' or ``near-zero'' if $|\sigma|\le\epsilon$. Moreover, given an integer $s\le n$, for a vector $\hat{x}\in\mathbb{R}^n$ having support size $\ge s$, we say that $\hat{x}$ ``almost'' has support size $s$ if
\begin{equation}
  |\supp{\max(\mathbf{0},|\hat{x}|-\epsilon\mathbf{1})}|=s. \label{eq:amb}
\end{equation}
\red{Notationwise, the inner $|\cdot|$ operator denotes absolute value, while the outer denotes set cardinality. Moreover, operators are applied vectorially}: $\mathbf{0}$ and $\mathbf{1}$ are the all-zero and all-one vectors, and we extend the absolute value and maximum operators $|\cdot|,\max$ to apply to vectors componentwise. So a vector almost has support size $s$ if there are exactly\footnote{Some results in \cite{damelin} relax \textit{exactly} to \textit{at most}. This relaxation is, however, motivated and formally explained.} $s$ components with absolute value larger than $\epsilon$.

\begin{thm}
  \label{thm:sparse}
  Given $\hat{x}\in\mathbb{R}^n$, a scalar $\delta\in(0,1/(2+2\sqrt{2}))$, \red{and a $m\times n$ matrix $A$ sampled componentwise from a normal distribution with mean $0$ and standard deviation $1$ (denoted $A\sim\mathsf{N}(0,1)^{mn}$)}, there exist two constants $c_1,c_2>0$ depending on $\delta$ such that, under the following assumptions:
  \begin{enumerate}
      \setlength{\parskip}{-0.2em}
  \item $\hat{x}$ almost has support size $s$;
  \item $m\ge \frac{1}{c_1}s\ln(\frac{n}{s})$;\label{thm:mLB}
  \item $P^1(A,A\hat{x})$ has a minimum $x^\ast$ with $|\supp{x^\ast}|=s$,\label{thm:uniq}
  \end{enumerate}
  we have
  \begin{equation}
    \Prob{\|x^\ast-\hat{x}\|_1\le 2\frac{1+2\delta(\sqrt{2}-1)}{1+2\delta(1-\sqrt{s})}\|\tilde{x}-\hat{x}\|_1\le (n-s)\epsilon}\ge 1-e^{-c_2m}. \label{goodapprox}
  \end{equation}
  \red{where $\tilde{x}$ is obtained by zeroing the $n-s$ smallest (in absolute value) components of $\hat{x}$}.
\end{thm}
This theorem is a convincing example of my accusations of vagueness. (a) It is based on some ``constants $c_1,c_2$ depending on $\delta$'' (but we are not told how to compute them). (b) It proves a probabilistic statement which approaches $1$ rapidly as $m$ increases, but we may not be able to control $m$ in practice: after all if $m=n-d$ is large, it means that the transmitted message is much larger than the original one: what if the communication channel is both noisy and costly? (c) It expects that a human reader would interpret the almost unparseable statement in Eq.~\eqref{goodapprox} as ``$x^\ast$ is a good approximation of $\hat{x}$'': what if it is not quite as good an approximation as we need? We can certainly fine-tune $\delta$, but probably at the expense of the constants $c_1,c_2$.

In fact, my understanding of these results is that they explain phenomena that were already known computationally within an engineering community. They provide a qualitative account, more than a workable recipe. These theories work inconclusively often, with unforeseeable exceptions that are possibly curable at the expense of increasing sizes to unwieldy millions, billions or more. Nonetheless, they have the merit of addressing foundational doubts that the engineering communities must have had. ``Goodness gracious, I'm designing a communication protocol for a new aircraft, and I'm using a technique which has always worked, but no-one knows why. I'm going to hell, aren't I?'' Such engineers will sleep soundly, from now on --- at least with arbitrarily high probability depending on an unknowable constant.

Proving Thm.~\ref{thm:sparse} requires quite a bit of time and patience, at least in the treatment of \cite{damelin}. The strategy is as follows:
\begin{enumerate}
    \setlength{\parskip}{-0.2em}
\item If $A$ has a certain complicated property called ``null space property'' (NSP), Eq.~\eqref{goodapprox} follows.
\item If $A$ has a certain other somewhat less complicated property called ``restricted isometry property'' (RIP), then $A$ also has the NSP.
\item If $A\sim\mathsf{N}(0,1)^{mn}$, then $A$ has the RIP.
\end{enumerate}
Note that this proof is a chain of \red{implications}. We end up proving that matrices sampled componentwise from a normal distribution are good, but this does not mean they are the only ones. Moreover, the lower bound for $m$ given in \red{Assumption \ref{thm:mLB} in Thm.~\ref{thm:sparse}} is valid, but it may not be tight. In proving Thm.~\ref{thm:sparse}, by the way, we also learn that the NSP guarantees that the minimum $x^\ast$ of $P^1(A,A\hat{x})$ having support size $s$ is unique.

Again, all this comes mostly from \cite{damelin}. The theoretical treatment of Thm.~\ref{thm:sparse} in \cite{matousek-lp} is more compact than that of \cite{damelin}, and has no probabilistic statements of the likes of Eq.~\eqref{goodapprox}, but it actually proves a much weaker and less general result than Thm.~\ref{thm:sparse} --- only enough to justify the LP with the application at hand. I find that the treatment in \cite[\S 5.5]{moitra} strikes a better balance between vagueness, clarity, and generality.

To go back to our original issue, we now know that solving $P^1(A,b)$ gives a good approximation to the sparsest solution of $Ax=b$ with high probability, as long as the number of rows $m$ of $A$ is large enough: we have a lower bound for $m$, but it is not tight.

\red{In case you are still wondering about the one-pixel camera, let me provide some closure.} I did find a few mentions of a ``one-pixel camera'' connected with the names usually attributed to Thm.~\ref{thm:sparse} (namely: ``compressed sensing'', ``compressive sampling'', ``mathematics of sparsity''). I would not be surprised if the physical existence of this object were just as vague as many of the theoretical statements that no doubt prompted my colleague to mention it to me.

\section{Finding hay in a haystack}
\label{s:randproj}
I want to go back to Eq.~\eqref{jll1}, and replace that nasty ``$\approx$'' sign by something more precise. Johnson and Lindenstrauss worked in a period where statements such as Eq.~\eqref{goodapprox} were less common. Their lemma is existential rather than probabilistic (the method of proof is probabilistic \cite{alon}, however). Again with minor notation changes, we have:
\begin{lem}
  \label{lem:jll}
  For each $\varepsilon\in(0,1)$ there is a constant $c=c(\varepsilon)>0$ (not depending on $n$) so that if $X\subset\mathbb{R}^n$ with $|X|=n$, then there is a mapping $f:X\to\mathbb{R}^k$ (where $k=\lceil c(\varepsilon)\ln n\rceil$) which satisfies
  \begin{equation}
    \forall x,y\in X\quad(1-\varepsilon)\|x-y\|_2\le\|f(x)-f(y)\|_2\le(1+\varepsilon)\|x-y\|_2. \label{jll2}
  \end{equation}
\end{lem}
Any probabilistic proof argues a positive probability of the existence of a certain entity, whence it infers that the entity must exist. The entity referred to in Lemma \ref{lem:jll} is the function $f$. The proof shows that $f$ is a random orthogonal projection to a lower-dimensional subspace --- of course it does not exclude that $f$ could take other forms.

Other than that, the proof is barely understandable. It contains sentences an old-school mathematician would never want to read in a proof, e.g.~``We plugged in the exact constant of $\sqrt{2}$ in Khintchine's inequality \red{\cite{khintchine}}, but of course any constant would serve as as well'', which the part of me who was educated on Kenneth Kunen's \textit{Set Theory} \cite{kunen} considers as inacceptably vague\footnote{This should not be misconstrued into a feeling that I lack respect for this result or either of its authors! I am only lamenting my lack of imagination when reading creative proofs.}. If the constants do not matter, how can one ever hope to use the result in practice, where computation requires that all constant symbols must take a definite value? Or are we to understand that we are facing a computation which will give the same result no matter what the value of the constant is? But that would simply imply that the statement is independent of such a constant, which should therefore be ignored: so why isn't it? And so on with the nagging doubts.

At the end of the proof we learn that
\begin{equation}
  c(\varepsilon)\ge \frac{100}{\varepsilon^2}. \label{jll:const}
\end{equation}
Now suppose you had a set $X$ of $1000$ vectors in $\mathbb{R}^{1000}$, and you wanted to lose some dimensions, obtaining a corresponding set $f(X)$ such that all pairwise distances are preserved up to $1\%$ accuracy. Then $\varepsilon=0.01$, which means that the projected dimension $k=\lceil 100/(0.01)^2\ln(1000)\rceil$, which has value \red{$6,907,756$}: a remarkable increase from $n=1000$. Of course, if you had a set of one billion vectors in one billion dimensions, $k$ would be just over one fiftieth of that, at $20,723,266$. From one trillion, the figure for $k$ would be only $28$ millions. When the logarithm in the growth order starts kicking in, we are obviously doing fine. The fact that the minimum value for $k$ starts at six millions (in the above example), however, makes one wonder about the practical usefulness. The issue, here, resides in the ``inacceptably vague'' statement in the proof of Lemma \ref{lem:jll} to the effect that Johnson and Lindenstrauss chose a certain constant (namely $\sqrt{2}$) resulting in the number $100$ in Eq.~\eqref{jll:const}. Choosing something other than $\sqrt{2}$ would have perhaps yielded a value smaller than $100$ in Eq.~\eqref{jll:const}.

Later proofs (see e.g.~\cite{dasgupta}) contributed a few other fundamental concepts to the Johnson-Lindenstrauss Lemma (JLL): (i) the notation $k=\lceil C\varepsilon^{-2}\ln(n)\rceil$, which makes the existence of coefficient $C$ explicit, without providing any hint about its value; (ii) the fact that the dimension $m$ of the vectors in $X$ need not be equal to the number $n$ of such vectors; (iii) an estimation of the probability of $f$ to do what it's supposed to (a probability which tends to one with exponential speed, but where the exponential term involves a ``universal constant'' $c$ which we cannot really compute); (iv) a variety of simpler constructions for $f$. Nowadays, we take $f$ to be $k\times n$ matrices sampled from sub-Gaussian distributions, of which the Gaussian distribution is one. But there are also nicer sub-Gaussian distributions, for example the distribution which samples $1$ or $-1$ with probability $1/6$ and $0$ with probability $2/3$ \cite{achlioptas}, its subsequent variant where $\prob{1}=\prob{-1}=\alpha/2$ and $\prob{0}=1-\alpha$ \red{for some given $\alpha<1$} \cite{kane}, or the more comfortable sparse normal matrices \red{with given sparsity} described in \cite[\S 5.1]{rpqp}.\label{alphapage}

The literature on the JLL is by now vast. In modern treatments (e.g.~\cite{vershynin}), the JLL is presented as a statement that has clear analogies with Thm.~\ref{thm:sparse}.
\begin{thm}
  \label{thm:jll}
  Let $X$ be a set of $n$ points in $\mathbb{R}^m$, $\varepsilon>0$, and $C,c$ be two universal constants. Then there exist an integer $k\ge\frac{C}{\varepsilon^2}\ln(n)$, and a $k\times m$ matrix $T$ sampled componentwise from a sub-Gaussian distribution with zero mean and standard deviation $\frac{1}{\sqrt{k}}$, such that:
  \begin{equation}
    \Prob{\forall x,y\in X\ (1-\varepsilon)\|x-y\|_2\le\|Tx-Ty\|_2\le(1+\varepsilon)\|x-y\|_2}\ge 1-2\red{n(n+1)}e^{-c\varepsilon^2k}.\label{jll3}
  \end{equation}
\end{thm}
In the form of Thm.~\ref{thm:jll}, the result says that any random sub-Gaussian matrix $T$ is a good choice for the $f$ in the JLL (Lemma \ref{lem:jll}). The issue of determining the constant has been made a non-issue: the JLL today is viewed as a principle more than a recipe for computation. It is up to the computer programmer to find values for $C,c$ that verify Eq.~\eqref{jll3}. Some attempts have been made in this sense, see e.g.~\cite{venkatasubramanian}. Once you try sampling $T$ you soon realize that values in the range $[0.5,2]$ seem to be fine, which lowers the usefulness bar from ``over one million'' to ``between $5,000$ and $20,000$'', a much more acceptable threshold. If you want smaller thresholds, you can allow for higher errors, and raise that 1\% to 10\% or even 20\%. There are many more caveats in order to successfully use the JLL in practice (see e.g.~\cite[\S 7.3.1]{dgds} \red{and \cite{rp4lp}}), but there is also considerable flexibility.

The vagueness of these probabilistic results turns out to be a blessing, albeit in disguise. It only captures a wealth of gauges which you can calibrate for your own uses. The trade-off of this flexibility is that making theoretically valid choices for the universal constants is only possible if the specific mathematical properties of the problem at hand allow it (moreover, you have to construct the theoretical arguments, which takes time and skill). Otherwise, if you lack skill and time (like me), you can resort to informed guessing aided by empirical tests.

There is one further safeguard against poor choices for $C$, $c$, $k$, and $T$ in Thm.~\ref{thm:jll}. The JLL is based on the phenomenon of ``concentration of measure'' \cite{ledoux2005,barvinok_710}, which ensures that a certain sampled function has high probability of being close to its median value. As a consequence, the event of sampling a random projector $T$ that \textit{completely fails} (i.e.~it fails for most of all pairs $x,y$) to satisfy its expected property Eq.~\eqref{jll2} \red{(where $f(x)=Tx$)} is going to be extremely rare. It may be the case that Eq.~\eqref{jll2} fails to hold for a few pairs $x,y$, but it will hold for most pairs. This may or may not be an issue, depending on the problem at hand. But in most large-dimensional datasets, where the data are prone to be partly wrong and noisy, small failures of Eq.~\eqref{jll2} should not impair the usefulness of the dimensionally reduced vector set $TX$.

In summary, acceptable random projectors look more like hay (than needles) in a haystack.

\section{Epiphany and communion}
\label{s:idea}
I wrote in Sect.~\ref{s:1norm} that the bound in \red{Assumption \ref{thm:mLB} in Thm.~\ref{thm:sparse}} is not tight. In choosing the $n\times d$ encoding matrix $Q$, it would be great to \red{decrease $n$}, so that the resulting encoded vector $z=Qw$ is not excessively costly in terms of storage and bandwidth: after all, the noisy channel is likely to be costly too, especially if it is a two-way communication link between Earth and a deep-space probe, as in Matou\v{s}ek and G\"artner's example. But what sort of result can I invoke in order to reassure myself that I am following, if not a theoretically proven recipe, at least a theoretically solid principle?

Let us go back to the formulation of $P^1(A,b)$ in Eq.~\eqref{basispursuit}: the equation constraints $Ax=b$ can be equivalently written as $\sum_j A_jx_j=b$, where $j\le n$ indexes the column vectors $A_j$ of $A$, and $x_j$ is the $j$-th component of the decision variable vector $x$. We know that $X=\{A_1,\ldots,A_n,b\}\subset\mathbb{R}^m$. If we apply the JLL to $X$, we should obtain lower dimensional vectors $TX=\{TA_1,\ldots,TA_n,Tb\}\subset\mathbb{R}^k$, with $k= O(\varepsilon^{-2}\ln(n+1))$. If $m$ is large enough, we can hopefully concoct some universal constants, and choose an appropriate $\varepsilon$, such that $k\ll m$ (or at least $k<m$), and that the probability of Eq.~\eqref{jll2} is still high. But one question remains: will the corresponding LP $P^1(TA,Tb)$ yield approximately the same optimal value and optimal solutions as $P^1(A,b)$, simply because the pairwise distances between the vectors in $TX$ are a good approximation of the corresponding distances in $X$?

This question was answered in the affirmative (vaguely, as usual in these types of results) in \cite{jllmor,rp4lp}. This justifies, at least in principle, the application of the JLL to $P^1(A,b)$, and its replacement with the smaller-sized $P^1(TA,Tb)$. Forward-thinking readers might make the following objection. Suppose you are able to compute a tight lower bound for the number of rows $m$ of $A$, for which Thm.~\ref{thm:sparse} still held; again, you would be confronted with an LP, and again you could apply the JLL to its columns in order to make them shorter. Would this not negate the tightness of your bound for $m$ and create a paradox?

The answer to this objection may reside in the observation that varying $m$ in $P^1(A,b)$ induces a phase transition w.r.t.~the truth value of Eq.~\eqref{goodapprox} \cite{tropp}. When $m$ is at its empirical lower bound, decreasing $m$ by a few units will cause $P^1(A,b)$ to produce completely dense solutions arbitrarily far from the original sparse error vector, with sharply rising probability \cite[Fig.~1]{tropp}. I will venture to say that this is compatible with the vague nature of the JLL, full of unknown constants and guessed $\varepsilon$'s. The JLL works in ``areas'' with comfortably large neighbourhoods of variation tolerance. If you are on the edge of a phase transition, it may well fail to apply. Supposing I were able to compute the universal constants exactly in the case of compressed sensing (which I am not), and that I wanted to work with a tiny $\varepsilon$ in order to carefully control the risk of overstepping the phase transition threshold, the JLL would perhaps tell me that $k$ would need to be as large as $m$ or even higher (due to the huge constant $C/\varepsilon^2$, since the $\ln(n)$ term increases too slowly to wreak such havoc).

But I have other weapons to deploy: for example, the fact that messages are strings of characters rather than real vectors.

\section{Pardoning more errors}
\label{s:pardon}
We assume that all messages have been segmented into textual pieces of $d$ characters (including the spaces), that a redundancy factor $R>1$ is given, and that a certain $n\times d$ matrix $Q$ (with $n=Rd$) is known to both sender and recipient. Every party also knows that the communication line has an error rate $\Delta$. Let us now look at the encoding and decoding algorithm. 

Given a string $\mu$ such as ``\texttt{I am a string}'', we have a function called \textsf{string2bitlist} which takes $\mu$ as input and returns a binary vector $w\in\{0,1\}^{d}$ as output, with $d=8|\mu|$ This function works by expressing the ASCII codes of the string characters into a base-2 representation. We assume that each ASCII code is in the range $\{0,\ldots,255\}$, which means that each character takes 8 bits of binary storage to be encoded. The vector $w$ is then considered a real vector in $\mathbb{R}^{d}$. The sender computes $z=Qw\in\mathbb{R}^n$, which will be sent over the $\Delta$-noisy (and somewhat costly) communication line.

Note that \textsf{string2bitlist} has an inverse function \textsf{bitlist2string} which takes a vector $w\in\{0,1\}^d$ (where $d$ is an integer multiple of $8$), splits $w$ into $d/8$ contiguous sub-vectors in $\{0,1\}^8$, which are then expressed in base 10 and interpreted as ASCII codes. The output of \textsf{bitlist2string} is the corresponding string $\mu$.

We make a few assumptions specific to the receiving end: the recipient has computed an $m\times n$ matrix $A$ where $m=n-d$ and $AQ=0$. Moreover, the recipient has implemented two functions: $\lfloor\sigma\rceil$, which rounds \red{the scalar} $\sigma$ to the nearest integer, and $\mathsf{cap}(\sigma,[L,U])=\min(\max(\sigma,L),U)$, which restricts $\sigma$ to lie in $[L,U]$. The decoding algorithm works as follows:
\begin{enumerate}
  \setlength{\parskip}{-0.2em}
\item receive $\bar{z}$ on the noisy communication line;
\item compute $b=A\bar{z}$;
\item solve $P^1(A,b)$ and obtain the optimal solution $x'$; \label{decstep3}
\item compute $z'=\bar{z}-x'$;
\item compute $w'=\mathsf{cap}(\lfloor((\transpose{Q}Q)^{-1}\transpose{Q} z')\rceil, [0,1])$; \label{decstep5}
\item return $\mu'=\mathsf{bitlist2string}(w')$.
\end{enumerate}
We also want to evaluate the difference between the decoded message $\mu'$ and the original message $\mu$. We therefore also compute $\mu_{\mbox{\sf\scriptsize err}}=\mathsf{dist}(\mu,\mu')$ for some ``distance function'' $\mathsf{dist}$ (for example the number of different characters). Note that, at Step \ref{decstep3}, we know that if the assumptions of Thm.~\ref{thm:sparse} hold, then the solution of $P^1(A,b)$ will be unique and a good approximation of the communication error vector \red{$\bar{x}=z-z'$} (which has density $\Delta$, and support size $\lfloor \Delta n\rceil$).

The rounding and capping operation at Step \ref{decstep5} is extremely forgiving to errors, in the sense that even if the number $m$ of rows in $A$ is smaller than required by Eq.~\eqref{thm:mLB}, resulting in $x'$ being dense, after the rounding and capping operation in Step \ref{decstep5} every component of $x'$ in $(-\infty, 0.5]$ will have been set to zero (with every other component being set to one). Since we started from binary vectors, no undue trick is being played. We exposed our sent vector $z$ to a lot more noise than was being logically necessary by embedding it in a real vector space; and now we are reducing that noise in the most common-sense possible way. This might allow us to guess the universal constants $C,c$ of Thm.~\ref{thm:jll} carelessly, and set $\varepsilon$ to a very tolerant 20\%, and still be able to decode messages correctly. Moreover, if we are sending plain text in natural language, which comes with its own redundancy, we might even be able to tolerate some errors and still understand the message.

\section{The outcome}
I set up my code with $R=4$, $\Delta=0.08$ (in accordance with \cite{matousek-lp}), so it could solve $P^1(A,b)$ as well as $P^1(TA,Tb)$ where $T$ is a sparse $(1,0,-1)$ sub-Gaussian random projector with density $\alpha=0.02$ (see page \pageref{alphapage}) and tolerance $\varepsilon=0.2$. I used the following algorithm in order to generate appropriately sized matrices $A,Q$ such that $AQ=0$:
\begin{enumerate}
  \setlength{\parskip}{-0.2em}
\item sample an $n\times n$ matrix $M$ componentwise from a uniform distribution (e.g.~on $[-1,1]$): $M$ will have full rank with probability 1 (since the probability that there should be fortuitous linear relations is proportional to the volume occupied by a subspace in the full space: zero);
\item find an eigenvector matrix of $M^\top M$: this provides an orthonormal basis of $\mathbb{R}^n$;
\item concatenate $d$ (column) eigenvectors to make $Q$;
\item stack $m=n-d$ (row) eigenvectors to make $A$.
\end{enumerate}
Clearly, $AQ=0$ by construction, since $(Q,\transpose{A})$ is an orthonormal $n\times n$ matrix. The corresponding implementation is fast, since there are robust and efficient implementations of pseudorandom sampling and spectral decomposition.

I then simulated the noisy message encoding and decoding for various initial segments of the beginning of Virgil's \textit{{\AE}neid}'s second canto:
{\it\begin{quote}
  Conticuere omnes intentique ora tenebant \\
  inde toro pater {\AE}neas sic orsus ab alto:\\
  ``Infandum, regina, iubes renovare dolorem: \\
  Troianas ut opes et lamentabile regnum \\
  eruerint Danai, qu{\ae}qu{\ae} ipse miserrima vidi, \\
  et quorum pars magna fui. Quis talia fando, \\
  Myrmidonum Dolopumve, aut duri miles Ulixi \\
  temperet a lacrymis? Et iam nox humida c{\ae}lo \\
  pr{\ae}cipitat, suadentque cadentia sidera somnos. \\
  Sed si tantus amor casus cognoscere nostros, \\
  et breviter Troi{\ae} supremum audire laborem, \\
  quamquam animus meminisse horret luctuque refugit,\\
  incipiam.''
\end{quote}}

\begin{figure}[!ht]
  \begin{center}
    \includegraphics[width=0.5\textwidth]{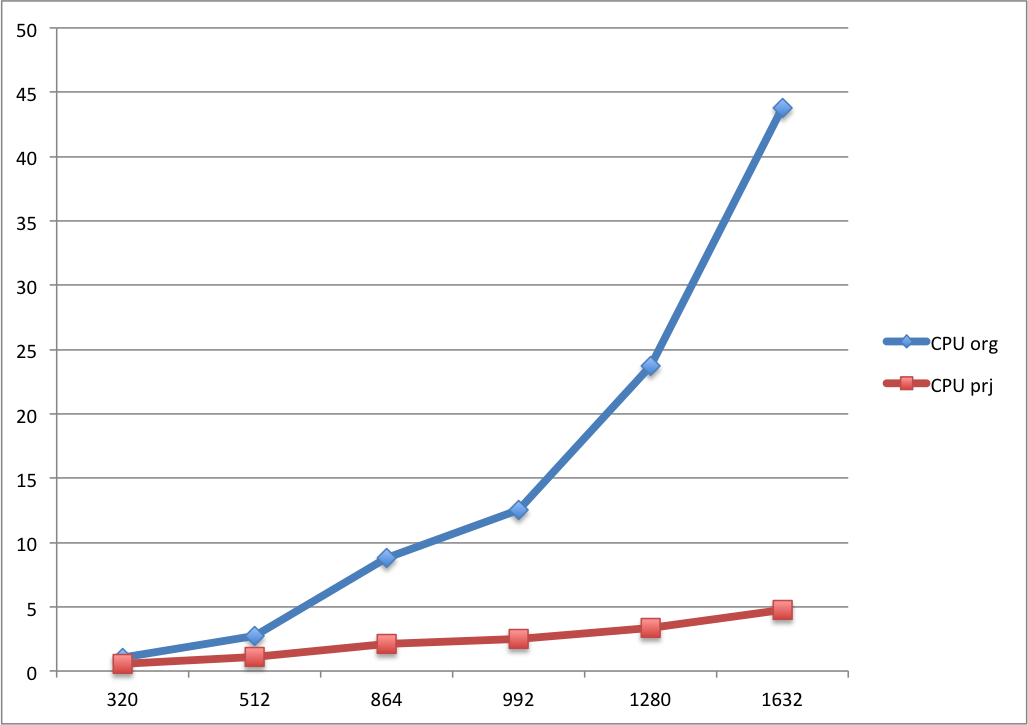}
  \end{center}
  \caption{CPU time for original and projected LP.}
  \label{f:res1}
\end{figure}

The computational results of these tests have been reported in Table \ref{t:res1}, where I reported the size ($d$ and $n$) \red{of each tested instance}, the error between original and decoded messages for the original LP $P^1(A,b)$ ($\mu^{\mathsf{org}}_{\mathsf{err}}$) and the projected LP $P^1(TA,Tb)$ ($\mu^{\mathsf{prj}}_{\mathsf{err}}$), with the corresponding CPU times. Evidently, the redundancy ratio $R=4$ for an error rate $\Delta=0.08$ (used in \cite{matousek-lp}) is overkill for computation, which is always successful even with the projected LP.
\begin{table}[!ht]
  \begin{center}
    \begin{tabular}{rr|rrrr}
      $d$ & $n$ & $\mu^{\mbox{\sf\scriptsize org}}_{\mbox{\sf\scriptsize err}}$ & $\mu^{\mbox{\sf\scriptsize prj}}_{\mbox{\sf\scriptsize err}}$ & {CPU${}^{\mbox{\sf\scriptsize org}}$} & {CPU${}^{\mbox{\sf\scriptsize prj}}$} \\ \hline
       80 &  320 & 0 & 0 & {1.05} & {\color{reasonablegreen}\textbf{0.56}} \\
      128 &  512 & 0 & 0 & {2.72} & {\color{reasonablegreen}\textbf{1.10}} \\
      216 &  864 & 0 & 0 & {8.83} & {\color{reasonablegreen}\textbf{2.12}} \\
      248 &  992 & 0 & 0 & {12.53} & {\color{reasonablegreen}\textbf{2.53}} \\
      320 & 1280 & 0 & 0 & {23.70} & {\color{reasonablegreen}\textbf{3.35}}\\
      408 & 1632 & 0 & 0 & {43.80} & {\color{reasonablegreen}\textbf{4.75}}
    \end{tabular}
  \end{center}
  \caption{Comparison of accuracy and CPU times taken to solve $P^1(A,b)$ and $P^1(TA,Tb)$. \red{Each line refers to a single run over the corresponding instance}.}
  \label{t:res1}
\end{table}
Since $P^1(TA,Tb)$ has fewer constraints than $P^1(A,b)$, it takes much less time to solve (see Fig.~\ref{f:res1}).

\section{Aiming for the impossible}
\label{s:impossible}
If you think about it, a redundancy $R=4$ for an error rate $\Delta=0.08$ appears excessive. I told myself that a good communication line should work with a redundancy $R$ that should be more or less $1+\Delta$. So then $n=(1+\Delta)d$. In particular, $n$ and $d$ would be very close (say $d=O(n)$), which would in turn make $A$ very short and fat: specifically, the number $m$ of rows of $A$ would be $n-d$, which would be very close to zero. On the other hand, Thm.~\ref{thm:sparse} requires $m$ to be ``large enough''. As we have seen before, ``almost zero'' can be ``large enough'' if all sizes involved are enormous. But I wanted this protocol to work with small sizes too, which means that, in practice, $m$ would need to be $O(n)$. So $m$ needs to be both ``close to zero'' and ``almost $n$'', clearly an impossible feat.

I decided to drop the first condition: let us make $m$ ``sizable'' with respect to $n$, say $O(n)$. Now the orthogonality condition $AQ=0$ implies that every one of the $O(n)$ rows of $A$ should be orthogonal to every one of the $O(n)$ columns of $Q$. This requires finding around $2n$ orthogonal vectors in $\mathbb{R}^n$, which is, again, impossible. Given the vagueness of the results leading to my idea, I thought it worthwhile to ask what one can obtain by relaxing ``orthogonality'' to ``near-orthogonality''.

As it happens, one can obtain quite a lot. Another way of seeing the JLL is that $\mathbb{R}^n$ can contain exponentially many almost orthogonal vectors. This follows from two corollaries of the JLL.
\begin{cor}
  Let $X$ be a set of $n$ points in $\mathbb{R}^m$, $\varepsilon>0$, and $C,c$ be two universal constants. Then there exist $k\ge\frac{C}{\varepsilon^2}\ln(n)$ and a $k\times m$ matrix $T$ sampled componentwise from a sub-Gaussian distribution with zero mean and standard deviation $\frac{1}{\sqrt{k}}$, such that:
  \begin{equation*}
    \Prob{\forall x,y\in X\ |\langle Tx,Ty\rangle - \langle x,y\rangle| \le \frac{\varepsilon}{2}(\|x\|_2^2+\|y\|_2^2)}\ge 1-4e^{-c(\varepsilon^2-\varepsilon^3)k}.
  \end{equation*}
\end{cor}
\begin{cor}
  Let $k\in\mathbb{N}$, $\varepsilon>0$, $C',c$ be two universal constants, and $n=\lceil e^{\frac{\varepsilon^2}{C'}k}\rceil$. Then
  \begin{equation*}
    \Prob{|\langle Te_i,Te_j\rangle| \le \varepsilon}\ge 1-4e^{-c(\varepsilon^2-\varepsilon^3)k}.
  \end{equation*}
\end{cor}
The first corollary says that if Euclidean distances are almost preserved, then vector angles are also almost preserved (it makes sense, as preserving distances leads to a congruence, which also preserves angles). The second corollary follows from the first corollary when it is applied to $X=\{e_1,\ldots,e_n\}=B_n$, the standard basis of $\mathbb{R}^n$, such that the constant $C'\le C$ is chosen in a way that yields $k=\lceil\frac{C'}{\varepsilon^2}\ln(n)\rceil$ in the first corollary: \red{then the second corollary says that there is sufficient space in $\mathbb{R}^k$ to host $O(e^k)$ almost orthogonal vectors.} Then the same will hold for $n$ by symbolic replacement.

The issue is now to see whether ``almost orthogonality'' is as good as orthogonality for the decoding process. After all, the last step of Eq.~\eqref{decode} states $(AQ)w + Ax = Ax$ because $AQ=0$. What if $\|AQ\|\le\varepsilon$ for some small $\varepsilon$ instead? A few tests with the almost orthogonal vectors that can be produced using $TB_{2n}$ with an $n\times 2n$ random projector $T$ led to disappointing results.

There is one further murky corner where I could try and scrape some more advantage for my impossible method: all I need, in order to decode $\bar{z}$ at least approximately, is that:
\begin{itemize}
 \item $A$ is an $m\times n$ matrix sampled componentwise from $\mathsf{N}(0,1)$ ($A$ will be used to solve $P^1(TA,Tb)$, where $n=\lfloor(1+\Delta)m\rceil$ and $T$ is a $k\times m$ random projector);
 \item $Q$ is an $n\times m$ matrix with $AQ\approx 0$: in other words, requiring that all the columns of $Q$ should be orthogonal is overkill.
\end{itemize}

I therefore formulated a LP in order to decide the components of $Q$. For every $\ell\le m$ let $q_\ell\in\mathbb{R}^n$ be the $\ell$-th column of $Q$. Then for each $\ell\le m$ I need to decide $q_\ell$ so that: (a) $Q=(q_\ell\;|\;\ell\le m)$ has full rank, and (b) $Aq_\ell\approx 0$. I attempted to enforce (a) by optimizing a random objective function, and (b) by imposing the linear constraints $Aq_\ell=0$. The following LP (Eq.~\ref{impossibleLP}) must be solved $m$ times \red{(each time with different random cost coefficients): the $\ell$-th solution $q_\ell$} must be stored in the corresponding column of $Q$ (for $\ell\le m$).
\begin{equation}
  \left.\begin{array}{rl}
    \max\limits_{q_\ell\in[-1,1]^n} & \sum\limits_{j\le n} \mathsf{Uniform}(-1,1)q_{\ell j}  \\
    \forall i\le m & \sum\limits_{j\le n} A_{ij} q_{\ell j} = 0. 
  \end{array}\right\}
  \label{impossibleLP}
\end{equation}
\red{Eq.~\eqref{impossibleLP} can be solved with any LP algorithm. In my experiments, I used the interior point method without crossing over to the simplex algorithm at the end, since I only needed feasible solutions.}

Let me add that solving these LPs in order to create $Q$ can take a considerable CPU time: definitely more than solving $P^1(A,b)$. But these matrices can be pre-computed for given sizes, and shared before communications occur. So I have not taken these CPU times into account in my tests. 

I was expecting Eq.~\eqref{impossibleLP} to lead to the inevitable \red{trivial solution $Q=0$ (which negates the requirement of $Q$ having full rank)}, and was prepared to increase the feasibility tolerances of the CPLEX solver I was using \cite{cplex128}. Instead, CPLEX provided me with a magnificent \red{(and impossible)} matrix $Q$ with full rank, and such that $AQ=0$ with precision $O(10^{-10})$. When I fine-tuned my rank evaluation code with sufficiently small numerical tolerances I discovered the \red{expected} rank deficiency, which was not acceptable since I would need $\transpose{Q}Q$ to be truly invertible for pseudoinverse purposes. So I switched $Q$ with $A$ in my process: I sampled $Q$ randomly, obtaining a full rank $Q$ with probability 1, and then applied Eq.~\eqref{impossibleLP} to the rows of $A$ instead. $A$ would not be a normally sampled matrix satisfying the assumption of Thm.~\eqref{thm:sparse}, but I hoped that the non-necessary nature of that theorem (which only lists sufficient conditions to its conclusion) would allow my $A$, \red{with full rank up to $O(10^{-10})$ precision and surely rank-deficient beyond, but} almost orthogonal to a random matrix, to do the trick nonetheless.

The results of my experiments are in Table \ref{t:res2}: on top of reporting the error measures for the accuracy of the two methods (original and projected LP) and the CPU times, I have also tried different values of $\Delta$. 
\begin{table}[!ht]
  \begin{center}
        \begin{tabular}{rrr|rrrr}
      $m$ & $n$ & $\Delta'$ & $\mu^{\mbox{\sf\scriptsize org}}_{\mbox{\sf\scriptsize err}}$ & $\mu^{\mbox{\sf\scriptsize prj}}_{\mbox{\sf\scriptsize err}}$ & {CPU${}^{\mbox{\sf\scriptsize org}}$} & {CPU${}^{\mbox{\sf\scriptsize prj}}$} \\ \hline
      328 & 426 & 0.3 & {182} & {\color{blue}\textbf{15}} & 2.45 & {\color{reasonablegreen} \textbf{1.87}} \\
      {\color{offwhite}328} & {\color{offwhite}426} & {\color{offwhite}0.3} & {154} & {\color{blue}\textbf{0}} & 2.20 & {\color{reasonablegreen} \textbf{1.49}} \\
      {328} & 459 & 0.4 & {\color{blue}\textbf{0}} & {1} & 4.47 & {\color{reasonablegreen} \textbf{2.45}} \\
      {\color{offwhite}328} & {\color{offwhite}459} & {\color{offwhite}0.4} & {\color{blue}\textbf{5}} & {17} & 2.86 & {\color{reasonablegreen} \textbf{1.46}} \\
      {328} & 492 & 0.5 & {60} & {\color{blue}\textbf{0}} & 4.53 & {\color{reasonablegreen} \textbf{1.18}} \\
      {\color{offwhite}328} & {\color{offwhite}492} & {\color{offwhite}0.5} & {34} & {\color{blue}\textbf{0}} & 5.38 & {\color{reasonablegreen} \textbf{1.18}} \\
      {328} & 590 & 0.8 & {14} & {\color{blue}\textbf{0}} & 8.30 & {\color{reasonablegreen} \textbf{1.41}} \\
      {\color{offwhite}328} & {\color{offwhite}590} & {\color{offwhite}0.8} & {107} & {\color{blue}\textbf{4}} & 6.76 & {\color{reasonablegreen} \textbf{1.43}} \\ 
      1896 & 2465 & 0.3 & {\color{blue}\textbf{0}} & {2} & 29.67 & {\color{reasonablegreen} \textbf{17.13}}
    \end{tabular}
  \end{center}
  \caption{Comparison of accuracy and CPU time taken to solve $P^1(A,b)$ and $P^1(TA,Tb)$ with the impossible choice of $A$. \red{Each instance was solved twice, and the results are reported for each of the two runs (black and grey entries)}.}
  \label{t:res2}
\end{table}

Table \ref{t:res2} tells an exhilarating and disturbing story. In most cases, solving the original LP $P^1(A,b)$ leads to Thm.~\ref{thm:sparse} failing, as can be ascertained from the high values of the error $\mu^{\mathsf{org}}_{\mathsf{err}}$. Solving the projected LP $P^1(TA,Tb)$, however, sometimes just does the trick. It is exhilarating to think that using one theorem wrongly is bad, but that using two of them wrongly is fine. It is even more satisfying to think that I am exploiting a redundancy ratio as small as the error ratio, and managing to reconstruct messages ``reasonably well''. On the other hand, it is disturbing to think that, in order to achieve these results, I have misused Thm.~\ref{thm:sparse}, Thm.~\ref{thm:jll}, basic linear algebra, and the LP machinery in ways that I really should not have done: \textit{eppur si muove}.


\section{Conclusion}
The lack of a formal (if vague) theoretical justification for the apparent success of my impossible method for decoding messages is part of the reason why this paper is written with an unusual style. Signal processing engineers, however, had used $P^1(A,b)$ instead of $P^0(A,b)$ for years, before a theory was constructed. So I am not despairing that one day someone will find a reason why this impossible method works. 

\section{Acknoledgements}
I am grateful to two anonymous referees for helping to improve this paper. 


\bibliographystyle{plain}
\bibliography{dr2}

\end{document}